\documentclass[10pt]{amsart}
\usepackage{amssymb}

\newcommand{\Q}{{\mathbb Q}}
\newcommand{\R}{{\mathbb R}}

\newcommand{\Sum}{\sum}

\newcommand{\F}{F({\mathbb R}^k \times I, n)}
\newcommand{\Lie}{{\mathcal L}}
\newcommand{\E}{Emb(I, {\mathbb R}^k \times I)}
\newcommand{\B}{{\mathcal B}}
\newcommand{\BT}{{\mathcal{BT}}}

\theoremstyle{plain}
\newtheorem{thm}{Theorem}[section]
\newtheorem{prop}[thm]{Proposition}
\newtheorem{lemma}[thm]{Lemma}
\newtheorem{cor}[thm]{Corollary}

\theoremstyle{definition}
\newtheorem{defn}{Definition}[section]

\theoremstyle{remark}

\begin{document}
\title{A one-dimensional embedding complex}
\date{October 27, 2000}
\author{Kevin P. Scannell} 
\address{Department of Mathematics and Computer Science\\
         Saint Louis University\\
         St. Louis, MO 63103} 
\email{scannell@slu.edu}
\thanks{The first named author was partially supported by NSF grant DMS-0072515}
\author{Dev P. Sinha}
\address{Department of Mathematics\\
         Brown University\\
         Providence, RI 02906}
\email{dps@math.brown.edu}

\begin{abstract}
\noindent 
We give the first explicit computations of 
rational homotopy groups of spaces of ``long knots'' 
in Euclidean spaces.  We define a spectral sequence
which converges to these rational homotopy groups
whose $E^1$ term is defined in terms of braid
Lie algebras.  For odd $k$ we establish a vanishing line for 
this spectral sequence, show the Euler characteristic
of the rows of this $E^1$ term is zero, and make calculations
of $E^2$ in a finite range.
\end{abstract}

\maketitle

\section{Introduction}

In this paper we introduce a spectral sequence which converges
to the rational homotopy groups of $\E$, for $k \geq 3$, which is
the space of embeddings 
of an interval in $\R^k \times I$
with fixed endpoints and tangent vectors at those endpoints 
(essentially, the space of long knots in $\R^{k+1}$).   
Our starting point is the work of \cite{Si00},
which defines such spectral sequences in terms of the topology
of configuration spaces. 
The paper \cite{Si00} in turn builds upon work of Goodwillie
and his collaborators \cite{GK00, GW99a, We99}, 
who have built a powerful theory for studying spaces of embeddings in general.

The rational homotopy groups of configuration spaces, which 
comprise the $E^1$ term, are
Lie algebras which are well-known (we call them ``braid
Lie algebras'').
Just as the study of cohomology of embedding spaces gives rise
to the study of graph cohomology, which has been studied extensively
\cite{BN95, Ko94, To00, Va92},
our complexes of braid Lie algebras are interesting new  
objects in quantum algebra.  
Similar complexes were described
by Kontsevich in his plenary talk \cite{Ko00}.

We start by reviewing the computation of the
rational homotopy groups of ordered configurations of points
in Euclidean space as a graded Lie algebra under Whitehead product, 
as well as some basics of free Lie algebras, which appear as 
subalgebras of these homotopy groups.  At that point,
we will have the necessary algebraic background to define
the chain complexes which are the rows of the
$E^1$ term of our spectral sequence.  It turns out that
through $E^2$, our spectral sequence for the homotopy
groups of $\E$ depends, up to regrading, only on the parity of $k$. 
We focus on odd $k$. 
We prove some fundamental facts about these 
complexes, such as the vanishing of their Euler characteristic.
We proceed to describe algorithms for computing 
the homology of these chain complexes,
and in the final section present the results of these computations in low
dimensions.
In some cases, the classes which arise in $E^2$ must survive, 
implying the existence of non-trivial spherical families
of embeddings.  Non-zero higher differentials are also possible.
We end with a brief discussion of the case of $k$ even, which 
pertains to the theory of finite-type knot invariants.

The second author would like to thank Tom Goodwillie for many 
helpful discussions and Ira Gessel for help in simplifying
the proof of Theorem~\ref{chi}.

\section{The Rational Homotopy Groups of Configuration Spaces} \label{compFn}

We remind the reader of the computations of the rational homotopy
groups of configuration spaces \cite{FN62}
and their Lie algebra structure under Whitehead product.
Throughout this paper, $\pi_\ast(X)$ will denote the homotopy groups
of $X$ tensored with the rational numbers.
Let $F(M, n)$ denote the space of ordered configurations
of $n$ distinct points in a manifold $M$.  We consider the projection
$\rho: F(M,n) \to F(M, n-1)$ defined by forgetting the last point in the
configuration, which is in fact a
fiber bundle whose fiber is 
$M \setminus \{ (n-1) \; {\rm points} \}$.  Let $\iota$
denote the inclusion of the fiber.  When $M =
\R^{k+1}$, the fibers are homotopy equivalent to $\bigvee_{n-1} S^k$,
and the projection map admits a section, by adding
a point (say in a fixed direction at a large distance) to a 
configuration of $n-1$ points.
This section leads to a splitting of the long exact sequence of a fibration
into split short exact sequences
$$0 \to \pi_i(\bigvee_{n-1} S^k) \to \pi_i (\F) \to 
\pi_i(F(\R^{k+1}, n-1)) \to 0.$$
By induction, we find that additively 
$$
\pi_i(\F) \cong \bigoplus_{j=1}^{n-1} \pi_i(\bigvee_j S^k).
$$

We now compute the structure of the rational homotopy groups $\pi_\ast(\F)$
as a Lie algebra under the Whitehead product.   

\begin{defn} \label{defB}
Let $\B^o_{n}$ (respectively $\B^e_n$) 
be the Lie algebra (super Lie algebra for $\B^e_n$)
generated over $\Q$ by classes $x_{ij}$ for 
$1 \leq i, j \leq n$ with relations
\begin{enumerate}
\item $x_{ij} = x_{ji}$ (respectively, $-x_{ji}$ for $\B^e_n$). \label{rel1} 
\item $x_{ii} = 0$ \label{rel2}
\item $[x_{ij}, x_{\ell m}] = 0 \; \; \text{if} \;\; \{ i, j \} \cap
\{\ell, m\} = \emptyset$ \label{rel3}
\item $[x_{ij}, x_{j\ell}] = [x_{j\ell}, x_{\ell i}] = [x_{\ell i}, x_{ij}]$.
                           \label{rel4}
\end{enumerate}
\end{defn}

We call $\B^o_n$ and $\B^e_n$ {\em braid Lie algebras}.

\begin{thm} \label{piconfig}
There is a Lie algebra isomorphism between 
$\pi_\ast(\F)$ and $\B^o_n$ if
$k$ is odd or $\B^e_n$ if $k$ is even.
\end{thm}

\begin{proof}

We first define classes which generate $\pi_\ast(\F)$ as a Lie algebra.
Pick a basepoint in $\F$, say with $z_i = (2i, 0, \ldots, 0)$ for definiteness.
There are $\binom{n}{2}$ generators of $\pi_k(\F)$,
corresponding to distinct pairs $\{i, j\} \subseteq \{1, \ldots, n\}$,
which we now realize geometrically.
We define $b_{ij} \in \pi_k(\F)$ as the class represented by the composite
of two maps.  First, we collapse $S^k$ onto $S^k \vee I$ by sending the 
``southern hemisphere'' of $S^k$ to $I$ through the height function.
Next, choose a path $\gamma_{ij}$ from $z_i$ to the point
$(2j - 1, 0, \ldots, 0)$ in the complement of the
other configuration points, and let $\iota_j$ denote the map which sends $S^k$
to the unit sphere about the point $z_j$.
To define $b_{ij}$ we compose the collapse map above with the map 
$S^k \vee I$ to $\F$ which sends
$t\in I$ to $\F$ as 
$(z_1, \ldots, z_{i-1}, \gamma_{ij}(t), z_{i+1}, \ldots, z_n )$
and $S^k$ to $\F$ as 
$v \mapsto (z_1, \ldots, z_{i-1}, \iota_j(v), z_{i+1}, \ldots, z_n )$.  

To see inductively that these classes are generators of
$\pi_k(\F)$, we simply
note that $b_{in}$ is equal to the image under $\iota_\ast$
of the generator of $\pi_k (\bigvee_{n-1} S^k)$ 
defined by the inclusion of the 
$i$th wedge factor.   

It is simple to check that these $b_{ij}$ satisfy the relations for
$x_{ij}$ in the definition of $\B^o_n$.
Note that from the usual graded commutativity
of the Whitehead product, brackets in $b_{ij}$ anti-commute when 
$k$ is odd and commute when $k$ is even.
Note that $b_{ij} = (-1)^{k+1} b_{ji}$
and $b_{ii} = 0$ so that relations (\ref{rel1}) and (\ref{rel2}) are
satisfied.

We next verify 
that the $b_{ij}$ satisfy relation (\ref{rel3}).  
Recall that if $\{f\}$ and $\{g\}$ are elements of $\pi_k(X)$
then $[\{f\}, \{g\}] = 0$ if and only if
$f \vee g \colon S^k \vee S^k \to X$ extends to $S^k \times S^k$.
If $\{ i, j \} \cap \{\ell, m \} = \emptyset$, the map
$b_{ij} \vee b_{\ell m}$ may be so extended by sending 
$$
v \times w \mapsto
(p_1, \ldots, p_{i-1}, \iota_{ij}(v), p_{i+1}, \ldots, p_{\ell-1},
\iota_{\ell m}(w), p_{\ell+1}, \ldots),
$$
where $\iota_{ij}$ is the composite of
the collapse map of
$S^k$ onto $S^k \vee I$ with $\iota_j \vee \gamma_{ij}$.
Informally we say that $z_i$ can travel around $z_j$ and $z_\ell$ 
can travel around
$z_m$ without having their paths (the images of $S^k$ under the projection
onto the $i$th and $\ell$th coordinates) intersect.

Next, we verify that the $b_{ij}$ satisfy relation (\ref{rel4}). 
Equivalently, we claim that $[b_{j \ell}, b_{ij} + b_{i\ell}] = 0$.
Informally, we say that $b_{ij} + b_{i\ell}$ 
is represented by a map in which $z_i$
travels around $z_j$ and $z_\ell$ but no other points in the configuration, and 
this may happen simultaneously as $z_j$ travels around $z_\ell$, giving an
extension of $(b_{ij} + b_{i\ell}) \vee b_{j\ell}$ to $S^k \times S^k$
similar to the one given for $b_{ij} \vee b_{\ell m}$.

We claim that relations (\ref{rel1}) through (\ref{rel4}) are a 
complete set of relations for $\pi_\ast(\F)$.  
This follows from the fact that these relations may be
used to reduce to an additive basis of Lie algebra monomials of the form 
$[ \cdots [ b_{im}, b_{jm}] \cdots b_{\ell m} \cdots ]$, where 
$i, j, \ell < m \leq n$.
We exhibit this claim algorithmically when discussing the
computations in \S \ref{S:Algs}; see in particular Algorithm 5.2.

\end{proof}

The fiber sequence above leads us to identify some subalgebras
of $\pi_\ast(\F)$ which are free Lie algebras.
Tensored with the rationals, the homotopy groups of wedges of spheres,
$\pi_{\ast+1}(\bigvee_j S^k)$ for $k>1$,
are well known \cite{Hi55, Wh78} to form
free Lie algebra under Whitehead product, with $j$ generators in degree
$\ast + 1= k $.  
Since the inclusion map $\iota : \bigvee_{n-1} S^k \to \F$ is injective on
homotopy, and by naturality of Whitehead 
products, the image of these homotopy groups under $\iota_\ast$ in
$\pi_\ast(\F)$ is a free Lie algebra which is generated by the 
classes $b_{in}$.

In the development of the spectral sequence 
we will in fact need the rational homotopy groups of 
$FT(k,n) = F(\R^k \times I, n) \times (S^k)^n$.   We call $FT(k,n)$ the
space of tangential configurations, thinking of the points in $S^k$ as
unit tangent vectors at points of a configuration. Recall that the homotopy
groups of a product of spaces is a direct sum of their homotopy groups,
and all Whitehead products between these summands are zero.  Let
$\Lambda^o$ (respectively $\Lambda^e$) be the free Lie algebra (respectively
super Lie algebra) on one generator.  Let $\BT^o_n$ denote 
$\B^o_n \oplus n\Lambda^o$ and similarly for $\BT^e_n$.

\begin{cor}\label{C:comp}
There is a Lie algebra isomorphism between $\pi_{\ast+1}(FT(k,n))$ and
$\BT^o_n$ if $k$ is odd or $\BT^e_n$ if $k$ is even.
\end{cor}
These isomorphisms respect the gradings involved.  We may grade $\BT^o_n$
according to the number of generators appearing in a bracket.
The $d$th graded summand of $\BT^o_n$ coincides with 
$\pi_{d(k-1) + 1}(FT(k,n))$.

\section{Free Lie algebras}


Let $\Lie(A)$ denote the free Lie algebra over $\Q$ on a set $A$ of 
symbols.  For our explicit computations, we must choose 
an additive basis for $\Lie(A)$.  Natural labels for elements
of free Lie algebras can be obtained from rooted,
planar binary trees (hereafter, referred to as simply a {\em trees}) 
with leaves labeled by 
elements of $A$.  Such a tree prescribes a bracketing of the
elements which label the leaves.
The number of leaves is the {\em degree} of the tree.
Trees with a root but no branches (degree one) are identified with the set
of symbols $A$.
When the context is clear, we will identify trees with
the free Lie algebra elements they produce.   The obvious product 
of two trees $x$ and $y$ (a tree with a new root, left subtree $x$,
and right subtree $y$) corresponds to the product in the Lie algebra
and will therefore also be denoted $[x,y]$.

A set $\mathcal{H}$ of trees 
is called a {\em Hall set} for $\Lie(A)$
\cite[\S 4.1]{Re93} if the following conditions hold:
\begin{enumerate}
\item $\mathcal{H}$ has a total order $\leq$
\item $A \subset \mathcal{H}$
\item If $h = [h_1,h_2] \in \mathcal{H}$ then 
$h_2 \in \mathcal{H}$ and $h < h_2$.
\item For any tree $h = [h_1,h_2]$ of degree at least two,
we have $h \in \mathcal{H}$ if and only 
if $h_1, h_2 \in \mathcal{H}$,
$h_1 < h_2$, and either $h_1 \in A$ or 
$h_1 = [x,y]$ with $h_2 \leq y$. 
\label{C:final}
\end{enumerate}

It is straightforward to show that
a Hall set forms an additive basis
of $\Lie(A)$ \cite[Thm. 4.9]{Re93} (cf. Algorithm 5.1 below). 
The basis elements 
comprising a fixed Hall set will 
be called {\em Hall trees}.
It is easy to see that (many) Hall sets exist
\cite[Prop. 4.1]{Re93}; for completeness, we give a 
quick description of an algorithm for creating one.
\vspace{2ex}
\par\noindent
{\bf Algorithm 3.1} (Generating a Hall set).
\textsf{
Given an ordered list of symbols $A$ and a positive 
integer $d$, this algorithm outputs a list $\mathcal{H}_d$ 
consisting of the elements of degree less than or equal to $d$
in a Hall set for $\mathcal{L}(A)$.
The list $\mathcal{H}_d$ will be sorted according 
to the total order on the Hall set.
}
\par
\begin{enumerate}
\item \textsf{
Set a counter $n = 1$.
Copy the list $A$ into $\mathcal{H}_d$.
}
\item \textsf{
If $n = d$, terminate the algorithm and output $\mathcal{H}_d$.
Otherwise, proceed.
}
\item \textsf{
Form all
products $[h_1,h_2]$ such that $h_1,h_2 \in \mathcal{H}_d$, $h_1 < h_2$,
the degree of $[h_1,h_2]$ is $n+1$, and condition (\ref{C:final}) is
satisfied in the definition of Hall set. 
The only choice to be made is
where to insert this new element in the ordering on $\mathcal{H}_d$;
for definiteness in performing the calculations
in \S \ref{S:calcs}, we insert $[h_1,h_2]$ into $\mathcal{H}_d$ 
as the immediate 
successor to $h_1$, thus $h_1 < [h_1,h_2] < h_2$ as required.
}
\item \textsf{
Increment $n$ and go to (2).
}
\end{enumerate}
\vspace{2ex}
\par\noindent
{\bf Example.} The output of Algorithm 3.1 with $A = \{a,b\}$ 
and $d = 5$ is the following list:
\begin{align*}
&a
\hspace{3ex}
[a,[[[a,b],b],b]]
\hspace{3ex}
[a,[a,[a,[a,b]]]]
\hspace{3ex}
[a,[a,[[a,b],b]]]
\hspace{3ex}
[a,[[a,b],b]] \\
&[a,[a,[a,b]]]
\hspace{3ex}
[a,[a,b]]
\hspace{3ex}
[[a,[a,b]],[a,b]]
\hspace{3ex}
[a,b]\\
&[[a,b],[[a,b],b]]
\hspace{3ex}
[[a,b],b]
\hspace{3ex}
[[[a,b],b],b]
\hspace{3ex}
[[[[a,b],b],b],b]
\hspace{3ex}
b
\end{align*}
\vspace{3ex}

The following result will be used in computing the Euler 
characteristic of the chain complexes which appear as the
rows of the 
$E^1$ term of our spectral sequence.

\begin{lemma} \cite[Cor. 4.14]{Re93} \label{L:count}
The number of Hall trees for $\mathcal{L}(A)$ of degree $d$ 
equals
$$
\frac{1}{d} \sum_{j | d} \mu(j) |A|^{d/j}
$$
where $\mu$ is the M\"obius function.
\end{lemma}

\section{The spectral sequence}

In this section we present an explicit realization of the 
spectral sequence introduced in \cite[\S 4]{Si00} which
converges to the rational homotopy groups of $\E$.
The spectral sequence arises from models of 
$\E$ which are reminiscent of cosimplicial spaces, but
whose combinatorics are based on Stasheff polytopes instead
of simplices.  The 
entries are (Fulton-MacPherson compactified versions of)
the ordered tangential configuration
spaces $FT(k,n) = F(\R^k \times I, n) \times (S^k)^n$

We can now describe the $E^1$ term of this spectral sequence.
We first describe an unreduced version (which will be 
denoted throughout by the addition of a tilde $\tilde{E}^1$), followed by
the reduced version (denoted simply $E^1$).  
Recall that one way to obtain the $E^1$ term
of the homotopy spectral sequence of a cosimplicial space
is by first passing to homotopy groups of the entries, which
if all entries are simply connected defines a cosimplicial
abelian group.  The $E^1$ term is then the chain complex
associated to this cosimplicial abelian group, which is bi-graded because
the homotopy groups themselves are graded.
We show in \cite{Si00} that even though our models
are based on Stasheff polyhedra,
applying homotopy groups to these models gives
rise to cosimplicial abelian groups.  Hence the $\tilde{E}^1$ term of
our spectral sequence is the chain complex of the
cosimplicial abelian group:
$$
\pi_\ast( FT(k, 0)) = pt. \overset{\Rightarrow}{\leftarrow} \pi_\ast(FT(k,1)) 
\overset{\Rrightarrow}{\Leftarrow} \pi_\ast(FT(k, 2)) \cdots$$

Here the coface maps $d^i_\ast$ are induced by maps $d^i$ on configuration
spaces (or rather their Fulton-MacPherson compactifications)
which are ``doubling'' the $i$th point in a tangential
configuration in the direction of the unit tangent vector determined by 
the $i$th factor of $S^k$, or 
if $i$ = $0$ or $n$ by adding a point to the configuration at 
$(\vec{0}, 0)$ or $(\vec{0}, 1) \in \R^k \times I$.  
The codegeneracy maps $s^i$ are defined by forgetting a 
point in the configuration.  

\begin{thm}[see \cite{Si00}] \label{T:dpsmain}
There is a second-quadrant spectral sequence whose $E^1$ term is given by
$\tilde{E}^1_{-p, q} = \pi_q(FT(k, p))$ and $d^1$ given 
by $\Sigma_i (-1)^i d^i_\ast$ which converges to $\pi_\ast(\E)$.
\end{thm} 

We now make the coface and codegeneracy maps algebraically explicit.
Recall from Corollary~\ref{C:comp} that
the rational homotopy groups of
$FT(k,n)$ are isomorphic to the Lie algebra $\BT^o_n$
(or $\BT^e_n$)
generated by 
classes $x_{ij}$ and $y_i$ for $1 \leq i,j \leq n$ and with relations
defined as in Definition~\ref{defB} and so that $[y_i, x_{j \ell}] = 0$
for all $i, j, \ell$ and $[y_i, y_j] = 0$ for $i \neq j$.

\begin{defn}
Define $\sigma^\ell(i)$ to be $i$ if $i<\ell$ and $i+1$ if $i> \ell$.
For $0 \leq \ell \leq n+1$ define $\partial^\ell \colon \BT^o_n \to \BT^o_{n+1}$
(respectively from $\BT^e_n$ to $\BT^e_{n+1}$) 
to be the Lie algebra homomorphism defined on generators as follows.
\[
\partial^\ell (x_{ij}) = 
\begin{cases}
x_{\sigma^\ell(i) \sigma^\ell(j)}               &\text{if $i,j \neq \ell$} \\
x_{i\sigma^\ell(j)} + x_{i+1,\sigma^\ell(j)}    &\text{if $i = \ell$}       

\end{cases}
\]

\[
\partial^\ell (y_i) =
\begin{cases}
y_{\sigma^\ell(i)}  &\text{if $i \neq \ell$} \\
x_{i,i+1} + y_i + y_{i+1}   &\text{if $i < j = \ell$}       
\end{cases}
\]

\end{defn}

\begin{defn}
For $1 \leq \ell \leq n$ define $\phi^\ell \colon \BT^o_n \to \BT^o_{n-1}$
(respectively from $\BT^e_n$ to $\BT^e_{n-1}$) 
to be the Lie algebra homomorphism defined on generators as follows.
\[
\phi^\ell(x_{ij})  =
\begin{cases}
x_{\sigma^\ell(i) \sigma^\ell(j)}    &\text{if $i,j \neq \ell$} \\
0            &\text{if $i$ or $j$ = $\ell$}
\end{cases}
\]
\[
\phi^\ell(y_{i})  =
\begin{cases}
y_{\sigma^\ell(i)}    &\text{if $i \neq \ell$} \\
0            &\text{if $i$ = $\ell$}
\end{cases}
\]

\end{defn}

The following proposition is immediate from the definitions of the
classes $b_{ij}$ and the maps $d^\ell$ and $s^\ell$. 

\begin{prop}
Under the isomorphisms of Corollary~\ref{C:comp} the homomorphisms
$d^\ell_\ast$ and $s^\ell_\ast$ coincide with $\partial^\ell$
and $\phi^\ell$ respectively.
\end{prop}

Making Theorem~\ref{T:dpsmain} algebraically explicit using
Corollary~\ref{C:comp} and the previous proposition leads us 
to the following spectral sequence whose $E^1$ term is defined
in terms of braid Lie algebras.

\begin{cor} \label{algvers}
There is a second-quadrant spectral sequence which converges to $\pi_\ast(\E)$
such that
$\tilde{E}^1_{-n, d(k-1)+1}$ is isomorphic to the $d$th graded
summand of $\BT^o_n$
(respectively $\BT^e_n$), 
$E^1_{-n,q} = 0$ when 
$q-1$ is not a multiple of $k-1$, and $d^1$ is given 
by $\Sigma_i (-1)^i \partial^i$.
\end{cor}

For the rest of the paper, we focus on the 
case in which $k$ is odd.

A useful reduction when studying cosimplicial abelian groups is
the replacement of the $n$th group by the intersection of the 
kernels of the codegeneracy maps.  Such a reduction does not change the
homology of the associated chain complex.
First note that the codegeneracy maps $\phi^\ell \colon \BT^o_n \to
\BT^o_{n-1}$ respect the direct sum decomposition 
$\BT^o_n = \B^o_n \oplus n\Lambda^o$.  Restricted to the $\Lambda^o$ factors,
the intersection of the kernel of the $\phi^\ell$ is zero unless $n$ 
is equal to one, in which case it is all of $\Lambda^o$. 
Restricted to the $\B^o_n$ factor,
the kernel of the codegeneracy map $\phi^n \colon \B^o_n \to
\B^o_{n-1}$ is the subalgebra generated by the classes
$x_{in}$, which is in fact a free Lie algebra (see the remarks following
the proof of Theorem~\ref{piconfig}).   We identify the
kernel of all of the $\phi^\ell$ as a submodule of this free
Lie algebra.

\begin{defn}
For $n>1$, let $M_{d,n}$ be the submodule of of the degree $d$ summand
of $\BT^o_n$ generated
by brackets of the classes $x_{in}$ such that 
each $i$ from $1$ to $n-1$ appears as an index.  Let $M_{1,1} = \BT^o_1$.
\end{defn}

\begin{thm} \label{reducedss}
There is a spectral sequence which converges to $\pi_\ast(\E)$
whose $E^1$ term is given by $E^1_{-n, d(k-1)+1} = M_{d,n}$
and whose $d^1$ is the restriction to this submodule of the
$d^1$ of Corollary~\ref{algvers}.
\end{thm}

Note that $M_{d,n} = 0$ for $d<n-1$, which leads to the following 
vanishing theorem.

\begin{cor} \label{vanish}
In the spectral sequence of Theorem~\ref{reducedss}, 
$E^1_{-p, q} = 0$ if $q<p(k-1) + 2 - k$.
\end{cor}

It is interesting to note that while the modules $M_{d, n}$ may be
defined purely in terms of the free Lie algebra (on $n-1$ generators), 
the boundary maps between them require extending the free Lie algebra
to a braid Lie algebra.  From the algebraic definition of $d^1$
it is not obvious that its restriction to $M_{d,n}$ maps to $M_{d,n+1}$.

Since computing the $E^2$ term amounts to computing 
the cohomology of the complexes $M_{d,\ast}$,
as a warmup we will compute the rank of $M_{d,n}$, which we 
denote $R(d,n)$, and will show that $\chi (M_{d,\ast}) = 0$.
Recall that the number of Hall trees of degree $d$ with $n$ symbols
is equal by Lemma~\ref{L:count} to  
$\frac{1}{d} \Sum_{j|d} \mu(j) n^{d/j}$.  We may produce a basis
of $M_{d,n}$ by first considering all brackets of degree $d$ and
throwing away ones in which fewer than $n$ elements appear.
We find that 
$$R(d,n) = \frac{1}{d} \Sum_{i=0}^n (-1)^i \binom{n}{i} 
        \Sum_{j|d} \mu(j) i^{d/j}.$$

We pause to define $S(d,n) = \Sum_{i=0}^n (-1)^i \binom{n}{i} i^d$,
which are essentially Stirling numbers.  
There is a combinatorial
interpretation of $S(d,n)$ as the number of surjections from a
$d$ element set onto an $n$ element set (to verify this, count all set
maps and subtract the non-surjections).  Note as well that the $S(d,n)$
have a generating function, as
$$\Sum_{m=0}^\infty S(m,n) \frac{x^m}{m!} = (e^x -1)^m.$$

Reordering the summations of $R(d,n)$ we find the following:

\begin{prop}  \label{P:count}
$R(d,n) = \frac{1}{d} \Sum_{j|d} \mu(j) S(d/j, n)$.
\end{prop}  

We may give
$R(d,n)$ a combinatorial interpretation in line with this equality
as the number of surjections of a $d$ element set to an $n$ element
set which are not invariant under any cyclic permutation of the $d$
element set, modulo cyclic permutations of the $d$ element set.
It would be interesting to find a bijection between such 
equivalence classes of surjections and a basis of $M_{d,n}$.
Such a combinatorial interpretation would be particularly
interesting for $M_{n,n}$ which, along with its $\Sigma_n$ action
by permuting the letters, is known as $Lie(n)$ and arises 
in the calculus of functors approach to homotopy theory
\cite{AM99}.

\begin{thm}\label{chi}
The Euler characteristic of $M_{d,\ast}$ is zero for $d > 2$.
\end{thm}

\begin{proof}
The Euler characteristic of the complex $M_{d,\ast}$ is by 
definition $\Sum_{\ell = 1}^d (-1)^\ell R(d,\ell)$, which  
after applying Proposition \ref{P:count}, 
reversing the order of summation, and ignoring zero terms,
is equal to
$$\frac{1}{d} \Sum_{j|d} \mu(j) \Sum_{\ell = 1}^{d/j} (-1)^\ell S(d/j, \ell).
$$
We claim that $\Sum_{\ell = 1}^m (-1)^\ell S(m, \ell) = (-1)^m$, which 
can be verified by computing the coefficient of $x^m/m!$ of 
$\Sum_{p = 0}^m (-1)^p (e^x - 1)^p$.  Hence the Euler characteristic 
is equal to $\frac{1}{d} \Sum_{j|d} \mu(j) (-1)^{d/j}$, which is
zero if $d > 2$.
\end{proof}

\section{Algorithms} \label{S:Algs}
In this section we provide a detailed description of the
methods used to compute the boundary operators in the 
complexes described above.   These algorithms can be 
performed by hand for the complexes of small degree $d$, 
but are best implemented on the modern electronic computer 
otherwise.

Because the product of two Hall trees is not 
necessarily a Hall tree, one must have an algorithm which takes 
an arbitrary tree representing a free Lie algebra element, 
and expresses it as a linear combination of Hall trees.
The proof that this algorithm terminates and produces
the desired result is contained in the proof of Theorem
4.9 in \cite{Re93}.

\vspace{2ex}
\par\noindent
{\bf Algorithm 5.1} (Hallification).
\textsf{
Given an integral linear combination of 
trees representing an element of $\Lie(A)$, 
this algorithm outputs a linear combination of Hall trees
representing the same element of $\Lie(A)$.
}
\par
\begin{enumerate}
\item \textsf{
If each tree appearing with a non-zero coefficient in the linear 
combination is Hall, terminate the algorithm and output the
linear combination.
Otherwise, choose $t$ to be the first non-Hall tree appearing
in the linear combination and proceed.
}
\item \textsf{
Find a subtree $s = [s_1,s_2]$ of $t$ which is
not Hall but whose children $s_1$ and $s_2$ are Hall.
This can be achieved by a simple recursion, noting that
the degree one trees (single letters) are Hall.
}
\item \textsf{
If $s_1 = s_2$, then remove $t$ from the linear combination
and go to step (1).
}
\item \textsf{
If $s_1 > s_2$, then switch $s_1$ and $s_2$ in $t$, 
multiply the coefficient of $t$ by
$-1$, and go to step (1).
}
\item \textsf{
We have $s_1 < s_2$.    
In this case, $s_1$ cannot be a single letter, or else
$s$ would be Hall.     So $s_1 = [x,y]$.  
We must have $y < s_2$ again using
the fact that $s$ is not Hall.   
Replace $t$ in the linear combination 
by the sum of two trees obtained by replacing
$s = [[x,y],s_2]$ by $[[x,s_2],y]$ and $[x,[y,s_2]]$
respectively and go to step (1).
}
\end{enumerate}
\vspace{2ex}

The following algorithm uses
the relations for $\B^o_{n}$ from Definition~\ref{defB}
and the Jacobi identity
to express elements
of $\B^o_{n}$ in a standardized form.
It will be used in the computation of the boundary
operator $\partial^n$ in Algorithm 5.3 below.

\begin{defn}
We say that a bracket in the classes $x_{ij}$ for $1 \leq i, j \leq n$
is {\em pure} if either all $x_{ij}$ which appear are of the form $x_{in}$
or none are of this form.
\end{defn}

\vspace{2ex}
\par\noindent
{\bf Algorithm 5.2} (Standard basis for $\B^o_n$).
\textsf{
Given an element $x$ of $\B^o_{n}$
expressed as a linear combination of brackets in the
$x_{ij}$, this algorithm computes a
linear combination of pure brackets also representing $x$.
}
\par
\begin{enumerate}
\item \textsf{
If each bracket appearing with a non-zero coefficient in the linear 
combination is pure, 
terminate the algorithm and output the
linear combination.
Otherwise, choose $t$ to be the first bracket in the linear combination 
which is not pure and proceed.
}
\item \textsf{
Find a smallest degree sub-bracket $s = [s_1,s_2]$ of $t$ which 
is not pure.
A simple recursion finds this sub-bracket.
}
\item \textsf{
If the degree of $s$ is two, go to step (4), otherwise go to step (7).
}
\item \textsf{
Since the degree of $s$ is two, we have $s_1 = x_{ij}$ and 
$s_2 = x_{\ell m}$ with either $j=n$ or $m=n$.
If $j=n$, go to step (5) and if $m=n$, go to step (6).
}
\item \textsf{
If $i=\ell$, then replace $t$ in the linear combination by
a new bracket obtained from $t$ by replacing $s = [x_{in}, x_{im}]$
with $[x_{mn}, x_{in}]$, using relation (4) in the
definition of $\B^o_n$. If $i=m$, then do the same thing,
replacing $s = [x_{mn}, x_{\ell m}]$ with 
$[x_{\ell n}, x_{mn}]$ by the same relation.   
In all other cases, remove $t$ from the linear combination 
(applying relation (3) in the definition of $\B^o_n$).
Start over at step (1).
}
\item \textsf{
If $i=\ell$, then replace $t$ in the linear combination by
a new bracket obtained from $t$ by replacing $s = [x_{ij}, x_{in}]$
with $[x_{in}, x_{jn}]$.  If $\ell=j$, then do the same thing,
replacing $s = [x_{ij}, x_{jn}]$ with 
$[x_{jn}, x_{in}]$.   In all other cases, remove
$t$ from the linear combination.  Start over at step (1).
}
\item \textsf{
If the degree of $s_1$ is greater than one, say $s_1 = [x,y]$,
we use the Jacobi identity
to replace $t$ in the linear combination by the sum of two  
brackets obtained from $t$ by replacing the sub-bracket $s = [[x,y],s_2]$ by
$[[s_2,y],x]$ and $[[x,s_2],y]$ respectively.  
If $s_1$ has degree one, then $s_2$ must have degree at least two,
say $s_2 = [x,y]$, and we do the same thing, replacing
$s = [s_1,[x,y]]$ by
$[x,[s_1,y]]$ and $[y,[x,s_1]]$ respectively.  
In either case, start over at step (1).
}
\end{enumerate}
\vspace{2ex}

A simple induction argument shows that this argument
terminates and produces the desired result.   
Namely, we associate to a bracket $t$ the
pair $(a,b)$ where $a$ is the number of generators 
$x_{ij}$ with $j < n$ appearing in $t$, 
and $b$ is the degree of the smallest impure sub-bracket
found in step (2).  We order such pairs 
lexicographically, with the minimum $(0,0)$
being achieved by pure brackets.  At every step, this
algorithm produces brackets whose associated pairs
are less than that of the original.  Steps (5) and (6),
corresponding to $b=2$, clearly reduce $a$.   Step (7)
leaves $a$ unchanged but reduces $b$, since
the sub-bracket $[x,y]$ of $s$ which is initially
pure becomes impure in all terms which occur after applying the
Jacobi identity.  Finally note that since $b$ in such an
associated pair is bounded by the degree of the bracket, 
there are only finitely many pairs less than a given one,
so the algorithm must terminate after a finite number of
recursive steps.

Note that the terms in the linear combination output by
Algorithm 5.2 which do not involve any $x_{in}$
can be run recursively through the algorithm 
as elements of $\B^o_{n-1}$,
yielding the standard form claimed in the proof
of Theorem \ref{piconfig}.

The final algorithm is the heart of the calculation; it
computes $\partial^\ell$ for $\ell = 0, \ldots, n$,
exploiting the fact that these maps are Lie algebra
homomorphisms.  Observe that $\partial^{n+1}$ is simply the 
natural inclusion of $\BT^o_n$ into $\BT^o_{n+1}$ 
and therefore requires no detailed description.

\vspace{2ex}
\par\noindent
{\bf Algorithm 5.3} (Boundary operator).
\textsf{
Given a basis element $t$ of $M_{d,n}$
(expressed as a Hall tree) and an integer $\ell$ between $0$ and $n$,
this algorithm computes $\partial^\ell(t)$ as a linear combination
of degree $d$ elements of $\B^o_{n+1}$ in the standard form
given by Algorithm 5.2.  
}
\par
\begin{enumerate}
\item \textsf{
If the degree of $t$ is greater than one, say $t = [t_1, t_2]$,
then recursively call Algorithm 5.3 to compute $\partial^\ell(t_1)$ and
$\partial^\ell(t_2)$.  Set $\partial^\ell(t) = [\partial^\ell(t_1), \partial^\ell(t_2)]$
and proceed to step (2).
If the degree of $t$ is one, go to step (4).
}
\item \textsf{
If $\ell = n$, then use Algorithm 5.2 to express the answer
$\partial^\ell(t)$ in standard form.  Proceed to step (3).
}
\item \textsf{
Use Algorithm 5.1 to express $\partial^\ell(t)$ in terms of Hall trees.
Terminate the algorithm and return $\partial^\ell(t)$.
}
\item \textsf{
If $\ell < n$, proceed to step (5), otherwise go to step (6).
}
\item \textsf{
Assume $t = x_{in}$.  
If $i < \ell$, set $\partial^\ell(t) = x_{i,n+1}$. 
If $i > \ell$, set $\partial^\ell(t) = x_{i+1,n+1}$. 
If $i = \ell$, set $\partial^\ell(t) = x_{i,n+1} + x_{i+1,n+1}$. 
Go to step (3).
}
\item \textsf{
Assume $t = x_{in}$. 
Set $\partial^\ell(t) = x_{in} + x_{i,n+1}$ and go to step (2).
}
\end{enumerate}
\vspace{2ex}

An example of this algorithm is worked out by hand
in the next section.

\section{Results} \label{S:calcs}
In this section we present some results of the computations
described in the previous section.   We will choose the gradings
to correspond to the case $k=3$, i.e. embeddings in $\R^4$.

First we note that in the degree one case, we have $E^1_{-1,3} = \Q$,
generated by $y_1$, $E^1_{-2, 3} = \Q$ generated by $x_{12}$, and $d^1$
is an isomorphism.
In degree two, the only non-zero entry is $E^1_{-3,5} = \Q$,
generated by $[x_{13}, x_{23}]$, implying $E^2_{-3,5} = \Q$.

We proceed by working out the first non-trivial boundary operator
$d^1 : E^1_{-3,7} \to E^1_{-4,7}$ by hand.
These spaces are by definition 
$M_{3,3}$ and $M_{3,4}$.   
Bases are obtained by creating, with Algorithm 3.1, Hall
bases for the free Lie algebra generated by $\{x_{13}, x_{23}\}$
(resp. $\{x_{14}, x_{24}, x_{34}\}$) 
and selecting the elements which have degree $3$ 
and such that all possible values of $i$ appear.
It turns out that each space is two-dimensional; the first is generated by
$[x_{13},[x_{13},x_{23}]]$ and $[[x_{13},x_{23}],x_{23}]$
and the second by 
$[x_{14},[x_{24},x_{34}]]$ and $[[x_{14},x_{34}],x_{24}]$.
Algorithm 5.3 is straightforward for $\ell \neq 3$;
in these cases we have:
\begin{align*}
&\partial^0[x_{13},[x_{13},x_{23}]] = [x_{24},[x_{24},x_{34}]]   \\
&\partial^0[[x_{13},x_{23}],x_{23}] = [[x_{24},x_{34}],x_{34}],
\end{align*}
while
\begin{align*}
\partial^1[x_{13},[x_{13},x_{23}]] &= [x_{14}+x_{24},[x_{14}+x_{24},x_{34}]]   \\
                                     &= [x_{14}+x_{24},[x_{14},x_{34}]+[x_{24},x_{34}]]   \\
                                     &= [x_{14},[x_{14},x_{34}]]+[x_{14},[x_{24},x_{34}]]+[x_{24},[x_{14},x_{34}]]+[x_{24},[x_{24},x_{34}]]   \\
                                     &= [x_{14},[x_{14},x_{34}]]+[x_{14},[x_{24},x_{34}]]-[[x_{14},x_{34}],x_{24}]+[x_{24},[x_{24},x_{34}]]   \\
\end{align*}
and
\begin{align*}
\partial^1[[x_{13},x_{23}],x_{23}] &= [[x_{14}+x_{24},x_{34}],x_{34}] \\
                                     &= [[x_{14},x_{34}]+[x_{24},x_{34}],x_{34}] \\
                                     &= [[x_{14},x_{34}],x_{34}]+[[x_{24},x_{34}],x_{34}], \\
\end{align*}
where the last line in the first case comes from an application 
of Algorithm 5.1 for the free Lie algebra over
$\{x_{14},x_{24},x_{34}\}$.
Similarly we have for $\ell = 2$
\begin{align*}
\partial^2[x_{13},[x_{13},x_{23}]] &= [x_{14},[x_{14},x_{24}]]+
                                        [x_{14},[x_{14},x_{34}]] \\
\partial^2[[x_{13},x_{23}],x_{23}] &= [[x_{14},x_{24}],x_{24}]+ 
                                        [[x_{14},x_{24}],x_{34}]+ 
                                        [[x_{14},x_{34}],x_{24}]+ 
                                        [[x_{14},x_{34}],x_{34}] \\
                                     &= [[x_{14},x_{24}],x_{24}]+ 
                                        [[x_{14},x_{34}],x_{34}]+
                                        2\ast[[x_{14},x_{34}],x_{24}]+ 
                                        [x_{14},[x_{24},x_{34}]], 
\end{align*}
where again the last line comes from Algorithm 5.1.
Finally, as noted above, $\partial^4$ is the natural inclusion:
\begin{align*}
&\partial^4[x_{13},[x_{13},x_{23}]] = [x_{13},[x_{13},x_{23}]]   \\
&\partial^4[[x_{13},x_{23}],x_{23}] = [[x_{13},x_{23}],x_{23}].
\end{align*}
The case $\ell = 3$ is much more computationally taxing,
as it requires the use of Algorithm 5.2:
\begin{align*}
\partial^3[x_{13},[x_{13},x_{23}]] &= [x_{13}+x_{14},[x_{13}+x_{14},x_{23}+x_{24}]] \\
                                   &= [x_{13},[x_{13},x_{23}]]+ 
                                      [x_{14},[x_{14},x_{24}]]+ 
                                      [x_{24},[x_{34},x_{14}]]+
                                      [x_{14},[x_{34},x_{24}]]  \hspace{4ex} \text{  (Alg. 5.2)} \\
                                   &= [x_{13},[x_{13},x_{23}]]+ 
                                      [x_{14},[x_{14},x_{24}]]+ 
                                      [[x_{14},x_{34}],x_{24}]-
                                      [x_{14},[x_{24},x_{34}]]  \hspace{4ex} \text{  (Alg. 5.1)} \\
\partial^3[[x_{13},x_{23}],x_{23}] &= [[x_{13}+x_{14},x_{23}+x_{24}],x_{23}+x_{24}] \\
                                   &= [[x_{13},x_{23}],x_{23}]+ 
                                      [[x_{14},x_{24}],x_{24}]+ 
                                      [[x_{14},x_{34}],x_{24}]+ 
                                      [[x_{24},x_{34}],x_{14}]  \hspace{4ex} \text{  (Alg. 5.2)} \\
                                   &= [[x_{13},x_{23}],x_{23}]+ 
                                      [[x_{14},x_{24}],x_{24}]+ 
                                      [[x_{14},x_{34}],x_{24}]- 
                                      [x_{14},[x_{24},x_{34}]]  \hspace{4ex} \text{  (Alg. 5.1)} \\
\end{align*}
Since $d^1 = \Sigma_i (-1)^i \partial^i$, we have from the above
calculations that
\begin{align*}
d^1[x_{13},[x_{13},x_{23}]] &= 0 \\
d^1[[x_{13},x_{23}],x_{23}] &= 2\ast[x_{14},[x_{24},x_{34}]] +
                                       [[x_{14},x_{34}],x_{24}]. 
\end{align*}
and so the matrix for the boundary operator with respect to our
chosen bases is given by
$\left(
\begin{smallmatrix}
0 & 2 \\
0 & 1
\end{smallmatrix}
\right)$.
We conclude that the boundary operator has rank one,
and so $E^2_{-3,7} \cong E^2_{-4,7} \cong \Q$.
Further (computer) calculations yield the $E^1$ and $E^2$ terms for $k$ odd 
given in Tables 1 and 2.

\begin{table}
\begin{center}
\begin{tabular}{c|c|c|c|c|c|c|c|}
$\Q^{120}$ & $\Q^{300}$ & $\Q^{260}$ & $\Q^{89}$ & $\Q^9$ &      &      & 13 \\
\hline
           &            &            &           &        &      &      & 12 \\
\hline
           & $\Q^{24}$  & $\Q^{48}$  & $\Q^{30}$ & $\Q^6$ &      &      & 11 \\
\hline
           &            &            &           &        &      &      & 10 \\
\hline
           &            & $\Q^6$     & $\Q^9$    & $\Q^3$ &      &      & 9  \\
\hline
           &            &            &           &        &      &      & 8  \\
\hline
           &            &            & $\Q^2$    & $\Q^2$ &      &      & 7  \\
\hline
           &            &            &           &        &      &      & 6  \\
\hline
           &            &            &           & $\Q$   &      &      & 5  \\
\hline
           &            &            &           &        &      &      & 4  \\
\hline
           &            &            &           &        & $\Q$ & $\Q$ & 3  \\
\hline
  -7 & -6  & -5 & -4 \hspace{1ex} & -3 \hspace{1.5ex}  & -2 \hspace{1ex} &  -1 \hspace{1ex} & \\
\hline
\end{tabular}
\vspace{2ex}
\caption{$E^1$ term for $k=3$}
\end{center}
\end{table}

\begin{table}
\begin{center}
\begin{tabular}{c|c|c|c|c|c|c|}
      $\Q$ & $\Q$   &      &      &      &      & 13           \\ \hline
           &        &      &      &      &      & 12           \\ \hline
           & $\Q^2$ & $\Q$ &      & $\Q$ &      & 11           \\ \hline
           &        &      &      &      &      & 10          \\ \hline
           &        &      &      &      &      & 9            \\ \hline
           &        &      &      &      &      & 8            \\ \hline
           &        &      & $\Q$ & $\Q$ &      & 7            \\ \hline
           &        &      &      &      &      & 6            \\ \hline
           &        &      &      & $\Q$ &      & 5            \\ \hline
           &        &      &      &      &      & 4            \\ \hline
           &        &      &      &      &      & 3            \\ \hline
      -7 \hspace{0.5pt}  & -6     & -5 \hspace{0.5pt}  & -4 \hspace{0.5pt}  & -3 \hspace{0.5pt}  & -2 \hspace{0.5pt}  &              \\ \hline
\end{tabular}
\vspace{2ex}
\caption{$E^2$ term for $k=3$}
\end{center}
\end{table}

These low-dimensional computations do not reveal any regular behavior.
Note that, as allowed because the Euler characteristic of the rows 
is zero, some rows vanish while most do not.  Note as well 
that there is no additional vanishing along the edge of the
vanishing line of Corollary~\ref{vanish}.  

All of the classes in Table 2 survive to $E^\infty$
except perhaps those in 
bidegrees $(-6, 13)$ and $(-3, 11)$ which could support a
$d^3$ differential.  

\begin{thm}
There are non-trivial classes in $\pi_n(Emb(I, \R^3 \times I))$ 
for $n = 2,3,4,5,6$.
\end{thm}
It would be interesting to find explicit spherical families of embeddings
which represent these classes.  One expects the evaluation map
$$\Delta^n \times \E \to F(\R^k \times I, n) \times (S^k)^n$$
to play a central role in relating these homotopy groups to
those of $F(\R^k \times I, n) \times (S^k)^n$ which appear in 
our spectral sequence.

We conclude with a brief description of 
some of the methods used to verify 
the computer calculations (beyond merely
computing examples by hand and comparing
with the computer output, which was done
extensively).  
Algorithm 3.1 was checked by an
independent function which verified that the 
generated trees were Hall, and checked the number 
of elements in the resulting Hall set against the 
dimension count given by Lemma~\ref{L:count}.
It was verified in the course of computing the $E^2$
term in Table 2 that $d^2 = 0$
for each of the chain complexes comprising $E^1$.
A similar mathematical fact which was not hard-coded 
into the application is that
the image of $M_{d,n}$ under $d^1$ lands in $M_{d,n+1}$
despite the fact that this is not the case for the 
individual homomorphisms $\partial^\ell$.
The ranks of the boundary operators were
verified using the linear algebra capabilities of a
symbolic mathematics package (Maple).
Finally, a nice check of the system as a whole was
provided by varying the algorithm for generating Hall
sets (noting the choices made in Algorithm 3.1) and
verifying that the ranks of all boundary operators 
remained unchanged.

\section{Further Work}

In further work \cite{SS01} we will investigate the case of $k$ even,
which includes the case of classical knots.
Though the spectral sequences of \cite{Si00} do not necessarily
converge, one can use those methods to produce knot invariants, 
which we show are of finite type.  In particular, an optimistic
view of rational homotopy theory predicts that the module
of classes
along the vanishing line of our spectral sequence (which for
$k = 2$ is the anti-diagonal) is isomorphic to the module of 
primitives in the 
Hopf algebra of finite type invariants \cite{BN95}.
To prove such a conjecture would involve relating the
combinatorics of braid Lie algebras to those of Feynman diagrams,
which could give a satisfactory explanation in terms of algebraic
topology of the appearance of Feynman diagrams in the study
of knots.

\bibliographystyle{amsplain}
\bibliography{data}

\providecommand{\bysame}{\leavevmode\hbox to3em{\hrulefill}\thinspace}
\begin{thebibliography}{10}

\bibitem{AM99}
G.~Arone and M.~E. Mahowald, \emph{The {G}oodwillie tower of the identity
  functor and the unstable periodic homotopy of spheres}, Invent. Math.
  \textbf{135} (1999), no.~3, 743--788.

\bibitem{BN95}
D.~Bar-Natan, \emph{On the {V}assiliev knot invariants}, Topology \textbf{34}
  (1995), no.~2, 423--472.

\bibitem{FN62}
E.~R. Fadell and L.~P. Neuwirth, \emph{Configuration spaces}, Math. Scand.
  \textbf{10} (1962), 111--118.

\bibitem{GK00}
T.~G. Goodwillie and J.~R. Klein, \emph{Excision statements for spaces of
  embeddings}, In preparation, 2000.

\bibitem{GW99a}
T.~G. Goodwillie and M.~S. Weiss, \emph{Embeddings from the point of view of
  immersion theory. {I}{I}}, Geom. Topol. \textbf{3} (1999), 103--118.

\bibitem{Hi55}
P.~J. Hilton, \emph{On the homotopy groups of the union of spheres}, J. London
  Math. Soc. (2) \textbf{30} (1955), 154--172.

\bibitem{Ko94}
M.~L. Kontsevich, \emph{Feynman diagrams and low-dimensional topology}, First
  {E}uropean {C}ongress of {M}athematics, {V}ol. {I}{I}, Progr. Math., vol.
  120, Birkh\"auser, Boston-Basel-Berlin, 1994, pp.~97--121.

\bibitem{Ko00}
\bysame, \emph{Operads of little discs in algebra and topology}, Lecture at the
  Mathematical Challenges Conference, UCLA, 2000.

\bibitem{Re93}
C.~Reutenauer, \emph{Free {L}ie algebras}, London Math. Soc. Monogr. (N.S.),
  vol.~7, Clarendon Press, Oxford, 1993.

\bibitem{SS01}
K.~P. Scannell and D.~P. Sinha, \emph{The calculus of embeddings and
  finite-type knot invariants}, In preparation, 2000.

\bibitem{Si00}
D.~P. Sinha, \emph{The topology of spaces of embeddings of the circle},
  Submitted, 2000.

\bibitem{To00}
V.~Tourtchine, \emph{Sur l'homologie des espaces de noeuds non-compacts},
  math.QA/0010017, 2000.

\bibitem{Va92}
V.~A. Vassiliev, \emph{Complements of discriminants of smooth maps: topology
  and applications}, Transl. Math. Monographs, vol.~98, Amer. Math. Soc.,
  Providence, 1992.

\bibitem{We99}
M.~S. Weiss, \emph{Embeddings from the point of view of immersion theory. {I}},
  Geom. Topol. \textbf{3} (1999), 67--101.

\bibitem{Wh78}
G.~W. Whitehead, \emph{Elements of homotopy theory}, Grad. Texts in Math.,
  vol.~61, Springer-Verlag, New York-Berlin-Heidelberg, 1978.

\end{thebibliography}
\end{document}